\documentclass[12pt, reqno]{amsart}
\usepackage{amsmath, amstext, amsbsy, amssymb, amscd}

\newtheorem{theorem}{Theorem}[section]
\newtheorem{lemma}[theorem]{Lemma}
\newtheorem{proposition}[theorem]{Proposition}

\theoremstyle{definition}

\theoremstyle{remark}
\newtheorem{remark}[theorem]{Remark}

\numberwithin{equation}{section}

\newcommand{\C}{ \mathbb C }
\newcommand{\ch}{\text{ch}}

\newcommand{\End}{{\rm End}}
\newcommand{\FFock}{{\mathcal F}}
\newcommand{\Hn}{{\mathbb H}_n}
\newcommand{\Hx}{{\mathbb H}_X}

\newcommand{\la}{\lambda}

\newcommand{\On}{{\mathcal O}^{[n]}}

\newcommand{\Tr}{ {\rm Tr} }
\newcommand{\vac}{|0\rangle}

\newcommand{\Xn}{ X^{[n]}}
\newcommand{\Z}{ \mathbb Z }

{\vskip-\lastskip\medskip
  \noindent
  {\em #1.}\enspace
  }%
{\qed\par\medskip
  }

\begin{document}
\title[Hilbert schemes, integrable hierarchies, and GW
theory] {Hilbert schemes, integrable hierarchies, \\and
Gromov-Witten theory}
\author[Wei-Ping Li]{Wei-Ping Li$^1$}
\address{Department of Mathematics, HKUST, Clear Water Bay, Kowloon, Hong
Kong } \email{mawpli@ust.hk}
\thanks{${}^1$Partially supported by the grant HKUST6170/99P}

\author[Zhenbo Qin]{Zhenbo Qin$^2$}
\address{Department of Mathematics, University of Missouri, Columbia, MO
65211, USA} \email{zq@math.missouri.edu}
\thanks{${}^2$Partially supported by an NSF grant}

\author[Weiqiang Wang]{Weiqiang Wang$^3$}
\address{Department of Mathematics, University of Virginia,
Charlottesville, VA 22904} \email{ww9c@virginia.edu}
\thanks{${}^3$Partially supported by an NSF grant}


\begin{abstract}
Various equivariant intersection numbers on Hilbert schemes of
points on the affine plane are computed, some of which are
organized into $\tau$-functions of $2$-Toda hierarchies. A
correspondence between the equivariant intersection on Hilbert
schemes and stationary Gromov-Witten theory is established.
\end{abstract}

\maketitle
\date{}

\section{Introduction}

The relations between the cohomology rings of Hilbert schemes of
points on algebraic surfaces and $\mathcal W$ algebras \cite{LQW}
through vertex operators have led to the following question (i.e.
Question~5 in Sect.~6 of \cite{QW}): what is the precise
connection between Hilbert schemes and integrable hierarchies? The
connection between the geometry of Hilbert schemes and vertex
operators is made through the so-called Chern character operators
\cite{LQW} (also cf. \cite{Lehn}). However the Chern character
operators are usually nilpotent for (cohomology) degree reasons.
This presents a serious difficulty for a sensible answer to the
above question in the framework of ordinary cohomology theory.

The main goal here is to initiate a direct link between {\em
equivariant} cohomology rings of Hilbert schemes of $n$ points
$\Xn$ on a quasi-projective surface $X$ and integrable
hierarchies, and to establish a correspondence with (stationary)
Gromov-Witten theory. In the present paper, we will treat the case
when $X$ is the affine plane, where some main idea is made clear.

A distinguished $T=\C^*$-action on $X=\C^2$ induces an action on
the Hilbert schemes $\Xn$ with finitely many fixed points (cf.
\cite{ES}). As explained by Vasserot \cite{Vas}, the Heisenberg
algebra construction in \cite{Na1} extends to the $T$-equivariant
cohomology of Hilbert schemes $\Xn$ (also cf. Nakajima
\cite{Na2}). In particular, all information of the equivariant
cohomology ring $H^*_T(\Xn)$ is encoded in a ring
$\Hn=H^{2n}_T(\Xn)$. The direct sum $\Hx = \oplus_{n\ge 0} \Hn$
becomes the bosonic Fock space of a Heisenberg algebra. The ring
$\Hn$ is further identified in \cite{Vas} with the class algebra
of the symmetric group $S_n$ (cf. \cite{LT, Wa} for the very
relevant study of these class algebras). Under such an
identification, we observe that the $k$-th equivariant Chern
characters of the tautological rank $n$ vector bundle over $\Xn$
correspond precisely to the $k$-th power-sum of Jucys-Murphy
elements, extending the earlier observation for ordinary Chern
characters in \cite{QW}. In particular, from the corresponding
result on the class algebra of symmetric group \cite{Wa}, the
Chern characters give rise to a set of ring generators for $\Hn$.
Therefore, one way to present the equivariant intersection theory
on $\Xn$ is to study the intersection numbers of these Chern
characters and their variants.

We further introduce the moduli spaces $\mathcal M(m,n)$, where
$m\in \Z, n \ge 0$. As varieties, $\mathcal M(m,n)$ is isomorphic
to $\Xn$. The equivariant cohomology ring of $\mathcal M(m,n)$
also naturally corresponds to a ring $\Hn^{(m)}$ which is
isomorphic to $\Hn$. We identify $\FFock = \oplus_{m,n} \Hn^{(m)}$
with the fermionic Fock space via the celebrated boson-fermion
correspondence. The pioneering work of Kyoto school (cf.
\cite{MJD} and the references therein) on connections among Fock
spaces, vertex operators, and soliton equations provides much
algebraic background for the geometric picture developed here.

The intersection numbers of the equivariant Chern characters in
$\mathcal M(m,n)$ can be organized into three types of generating
functions. The so-called {\em $N$-point function} which we can
compute has a simple relation with the $N$-point {\em
disconnected} series of stationary Gromov-Witten invariants of
$\mathbb P^1$. Note that $N$-point disconnected series is somewhat
more complicated than the $N$-point connected series which were
computed in \cite{OP}. We relate the second generating function
(called the {\em multi-point trace function}) in a simple and
precise form to the characters on the fermionic Fock space. The
latter has been computed by Bloch and Okounkov \cite{BO}, and it
also computes the stationary Gromov-Witten invariants of an
elliptic curve according to Okounkov and Pandharipande \cite{OP}.
This trace function is also intimately related to the trace
functions in the theory of vertex algebras first studied by Zhu
\cite{Zhu}. Yet another generating function which involves the
equivariant intersection numbers of the moduli spaces $\mathcal
M(m,n)$ is shown to be the $\tau$-functions for the $2$-Toda
hierarchies \cite{UT}. It is interesting to compare with \cite{OP}
where $\tau$ functions of $2$-Toda hierarchies were constructed
from Gromov-Witten invariants and Hurwitz numbers.

The notion of the {\em equivariant Chern character operator},
which gives rise to a master operator $\mathfrak H (z)$ acting on
$\FFock$, underlies the calculations of all the above generating
functions. In the study of the class algebras of symmetric groups,
a counterpart of this operator has also played a distinguished
role \cite{LT, Wa}. The operator $\mathfrak H (z)$ turns out to be
related to another operator $\varepsilon_0(z)$, which has played a
key role in the study of stationary Gromov-Witten theory
\cite{OP}, by the following formula:

\begin{eqnarray}
\mathfrak H (z) &=& \frac1{e^{z/2} -e^{-z/2}}
\left(\varepsilon_0(z) -\frac1{e^{z/2} -e^{-z/2}} \text{I}
\right). \label{eq:GWH}
\end{eqnarray}
The formula (\ref{eq:GWH}) defines the Gromov-Witten/Hilbert
correspondence. Combining with the Gromov-Witten/Hurwitz
correspondence in \cite{OP}, we also have a Hurwitz/Hilbert
correspondence. In fact, the same combinatorics of the symmetric
groups underlies these distinct geometric studies. In particular,
the use of Jucys-Murphy elements could also be used to clarify the
notion of completed cycles considered in \cite{OP}.

An important open question is to establish a direct geometric
connection behind the correspondence between (equivariant)
intersection theory of Hilbert schemes and (equivariant)
Gromov-Witten theory. It is possible that mirror symmetry and
connections with string theory may play a key role (cf. e.g.
\cite{Dij}). Extension of the results of this paper to other
quasi-projective surfaces including the total space of the
cotangent bundle $T^*\mathbb P^1$ will be presented elsewhere.

The plan of the paper is as follows. In Sect.~\ref{sec:setup}, we
set up the notations and review some known results. In
Sect.~\ref{sec:chern}, we formulate the Chern character operators
and introduce the key operators on the fermionic Fock space. In
Sect.~\ref{sec:Npoint}, we study the $N$-point functions of
intersection numbers on Hilbert schemes and the multi-point trace
functions. In Sect.~\ref{sec:tau}, we formulate the connections
with $\tau$-functions in $2$-Toda hierarchies.

{\bf Acknowledgments.} The project was initiated in \cite{QW} in
an effort to understand the connection between Hilbert schemes and
integrable hierarchies via vertex operators. Qin and Wang thank
Hong Kong UST for its warm hospitality and support. 
\section{The equivariant setup for Hilbert schemes}
\label{sec:setup}

\subsection{The torus action on the Hilbert schemes}

Let $T = \C^*$, and let $\theta$ be the $1$-dimensional standard
module of $T$ on which $s \in T$ acts as multiplication by $s$,
and let $t$ be the associated character. Then the representation
ring $\mathcal R(T)$ is isomorphic to $\mathbb Z[t, t^{-1}]$.

Let $X$ be an algebraic variety acted  by $T$. Let $H^*_T(X)$ and
$H^T_*(X)$ be the equivariant cohomology and the equivariant
homology with $\C$-coefficient respectively. Note that
$H^*_T(pt)=H^*(BT)=\mathbb C[t]$. Then there exist linear maps
\begin{eqnarray*}
\cup: \,\, H^i_T(X) \otimes H^j_T(X) \to H^{i+j}_T(X), \qquad
\cap: \,\, H^i_T(X) \otimes H^T_j(X) \to H_{i+j}^T(X).
\end{eqnarray*}
If $X$ is of pure dimension, then there exists a linear map $D:
H^i_T(X) \to H^T_i(X)$. If $X$ is smooth of pure dimension, $D$ is
an isomorphism. When $f: Y \to X$ is a $T$-equivariant and proper
morphism of varieties, we have a Gysin homomorphism $f_!: H_*^T(Y)
\to H_*^T(X)$ of equivariant homology. Moreover, when both $Y$ and
$X$ are smooth of pure dimension, we have the Gysin homomorphism
$D^{-1} \, f_! \, D: H^*_T(Y) \to H^*_T(X)$ of equivariant
cohomology, which will still be denoted by $f_!$.

Given an algebraic surface $X$, the Hilbert schemes of $n$ points
on $X$ is a nonsingular complex variety of dimension $2n$, cf.
\cite{Na1} for a general reference.

>From now on, let $X=\mathbb C^2$. The torus $T=\mathbb C^*$ acts
on the affine coordinate functions $w$ and $z$ of $X$ by $s(w,
z)=(s w, s^{-1}z)$. It induces an action on the Hilbert scheme
$X^{[n]}$ of $n$ points on $X$ with finitely many fixed points
parametrized by partitions of $n$ \cite{ES}. Let $\lambda$ be a
partition of $n$ and $\xi_{\lambda}$ be the fixed point on
$X^{[n]}$ corresponding to the partition $\lambda$. Let
$i_\lambda$ be the inclusion map $\xi_{\lambda}\to X^{[n]}$. Let
$1_{\xi_{\lambda}}\in H^0_T(\xi_{\lambda})$ be the unit, and thus
$[\xi_{\lambda}]=i_{\lambda!}(1_{\xi_{\lambda}})$ lies in
$H^{4n}_T(X^{[n]})$. Here and below, $[-]$ denotes the equivariant
fundamental cycle or its associated equivariant cohomology class.

Denote by $\mathbb C[t]'$ the localization of the ring $\mathbb
C[t]$ at the ideal $(t-1)$, and denote

$$\iota_n = \bigoplus\limits_{\lambda}
i_{\lambda}: (\Xn)^T \to \Xn.$$ We define $H^*_T((X^{[n]})^T)'
=H^*_T((X^{[n]})^T)\otimes_{\mathbb C[t]} \mathbb C[t]'$ and
define $H^*_T(X^{[n]} )'$ similarly. We denote

$$\iota_{n!} \colon H^*_T((X^{[n]})^T)'
\longrightarrow H^*_T(X^{[n]} )'$$
to be the induced Gysin map. By the localization theorem,
$\iota_{n!}$ is an isomorphism. The inverse $(\iota_{n!})^{-1}$ is
given by
\begin{eqnarray*}
\alpha \mapsto \left ( \frac{(i_{\lambda})^*(\alpha)}
{e_T(T_{\xi_{\lambda}})} \right )_{ \lambda }
\end{eqnarray*}
where $e_T(T_{\xi_{\lambda}})$ denotes the $T$-equivariant Euler
class of the tangent bundle of $\Xn$ over $\xi_\la$.

We define a bilinear pairing as in Vasserot \cite{Vas}

$$\langle -, - \rangle_n:
 H^*_T(X^{[n]})'\otimes_{\mathbb C[t]'}H^*_T(X^{[n]})'\to \mathbb
C[t]' $$ by

\begin{eqnarray*}
\langle \alpha, \beta\rangle_n =(-1)^n
p_{n!}(\iota_{n!})^{-1}(\alpha\cup\beta)
\end{eqnarray*}
where $p_n$ is the projection of the set $(X^{[n]})^T$ of
$T$-fixed points to a point.

\subsection{The ring $\Hn$}

For $k\ge n$, from the spectral sequence computation and
$H^{2k}(X^{[n]})=0$, we see that $H^{2k}_T(X^{[n]})=t^{k-n}\cup
H^{2n}_T(X^{[n]})$. Following Vasserot \cite{Vas}, we have an
induced product $\star$ on $\Hn \stackrel{\text{def}}{=}
H^{2n}_T(\Xn)$ such that

$$t^n \cup(x \star y)=x\cup y.$$
There is also an induced non-degenerate bilinear form $\langle
-,-\rangle_n: \Hn \otimes \Hn \rightarrow \C$. We define

$$\Hx = \bigoplus_{n=0}^\infty \Hn.$$

The bilinear forms $\langle -,-\rangle_n$ on $\Hn$ for all $n$
induce a bilinear form on $\Hx$, which is denoted by $\langle
-,-\rangle$. Given a linear operator $\mathfrak f \in \End (\Hx)$,
we denote by $\mathfrak f^*$ the adjoint operator with respect to
this bilinear form. The unit in $H_T^0(X^{[0]})$ will be denoted
by $\vac$.

For a fixed point $\xi_{\lambda}$ of $\Xn$, define $[\lambda] \in
\Hn $ as in \cite{Vas} by

\begin{eqnarray}  \label{rescaling}
 t^n \cup [\lambda]=(-1)^n h(\lambda)^{-1} [\xi_{\lambda}]
\end{eqnarray}
where $h(\la) =\prod_{\square \in \la} h_\square$ is the product
of the hook numbers associated to $\la$.

\begin{remark}  \label{rem:same}
Passing from $H^*_T(X^{[n]})'$ to $\Hn$ does not lose information
since we can recover the cup product and bilinear form on
$H^*_T(X^{[n]})'$ from those on $\Hn$. Thus, understanding the
ring $\Hn$ is the same as understanding the equivariant
intersection theory on $\Xn$.
\end{remark}

\subsection{The Heisenberg algebra}

Let $i > 0$. Denote by $Y$ the $T$-invariant subspace $0 \times \C
\subset X$. We define

$$Y_{n, i} = \{(\xi, \eta) \in X^{[n+i]} \times X^{[n]} \mid \eta
\subset \xi, \text{Supp}(I_\eta/I_\xi)=\{y\} \in Y \}.$$
Let $p_1$ and $p_2$ be the projections of $X^{[n+i]} \times
X^{[n]}$ to the two factors. We define the linear operator
$\mathfrak p_{-i} \in \text{End}(\Hx)$ by

\begin{eqnarray*}
\mathfrak p_{-i}(\alpha) = D^{-1}p_{1!}(p_2^*\alpha \cap [Y_{n,
i}]) \in \mathbb H_{n+i}
\end{eqnarray*}
for $\alpha \in H^{2n}_T(\Xn)$. Note that the restriction of $p_1$
to $Y_{n, i}$ is proper. We define $\mathfrak p_{i} \in
\text{End}(\Hx)$ to be the adjoint operator of $\mathfrak p_{-i}$.
Alternatively, letting $p_2'$ be the projection of $(X^{[n]})^T
\times X^{[n-i]}$ to $X^{[n-i]}$, we see that

\begin{eqnarray*}
\mathfrak p_{i}(\alpha) = (-1)^{i} \cdot D^{-1}p_{2!}'(\iota_n
\times \text{Id})_!^{-1} (p_1^*\alpha \cap [Y_{n-i, i}]) \in
\mathbb H_{n-i}
\end{eqnarray*}
for $\alpha \in H^{2n}_T(X^{[n]})$. Finally, we put $\mathfrak
p_0=0$. An argument parallel to the one in \cite{Na1} leads to the
following, cf. \cite{Vas}\footnote{There are sign typos in
\cite{Vas} on the formula for the annihilation operators and
Heisenberg commutation relations.}.

\begin{proposition}
The operators $\mathfrak p_n$, $n\in \Z$, acting on $\Hx =
\oplus_{n=0}^\infty \Hn$ satisfy the following Heisenberg
commutation relation:

$$[\mathfrak p_m, \mathfrak p_n ] = m\delta_{m,-n} \text{\rm I}.$$
Furthermore, the space $\Hx$ becomes the Fock space (i.e.
irreducible module) over the Heisenberg algebra with highest
weight vector $\vac$.
\end{proposition}

Given a partition $\la =(1^{m_1}2^{m_2} \ldots)$ of $n$, we denote
${\mathfrak z}_{\la} =\prod_{r \ge 1} r^{m_r} m_r!$, which is the
order of the centralizer of an element in $S_n$ of cycle type
$\la$. We then define

$$\mathfrak p_{-\la} = \frac1{{\mathfrak z}_\la}
\prod_{r \ge 1} {\mathfrak p_{-r}^{m_r}} \vac.$$
The $\mathfrak p_{-\la}$'s, as $\la$ runs over all partitions of
$n$, form a linear basis of $\Hn$. One has

$$ \langle \mathfrak p_{-\la}, \mathfrak p_{-\mu} \rangle
 = \frac1{{\mathfrak z}_\la}\delta_{\la,\mu}. $$
\section{The equivariant Chern character operator}
\label{sec:chern}

\subsection{The equivariant Chern characters}

Let $\mathcal Z_n =\{(\xi,x) \in \Xn \times X \mid x \in
\text{supp}(\xi) \}$ be the codimension 2 universal subscheme.
Denote by $\On$ the tautological rank $n$ vector bundle $\pi_{1*}
(\mathcal O_{\mathcal Z_n} \otimes \pi_2^* \mathcal O) =\pi_{1*}
(\mathcal O_{\mathcal Z_n})$ on $\Xn$ induced from the trivial
line bundle $\mathcal O$ on $X$, where $\pi_1, \pi_2$ denote the
projections of $\Xn \times X$ to the factors. Clearly $\On$ is
$T$-equivariant over $\Xn$.

The computation of the torus action on $\Xn$ (cf. e.g. \cite{Na2})
implies that

\begin{eqnarray} \label{f1}
\On |_{\xi_\lambda} =\bigoplus\limits_{\square\in
\lambda}\theta^{c_{\square}}
\end{eqnarray}
as a $T$-module, where $c_{\square}$ is the content of the box
$\square$ in the Young diagram associated to $\lambda$. Denote by
$\ch_{k,T}^{[n]}$ the $k$-th $T$-equivariant Chern character of
$\On$. In particular the zero-th Chern character $\ch_{0,T}^{[n]}$
equals the rank of the vector bundle $\On$, which is $n$. Then
$$\ch_{k,T}^{[n]}|_{\xi_{\lambda}}
=\frac{1}{k!}\sum\limits_{\square\in \lambda}( c_{\square}t)^k.$$

By the projection formula, we have in $H^*_T(X^{[n]})'$ that

\begin{eqnarray}
\ch_{k,T}^{[n]} \cup [\xi_{\lambda}] =\frac{1}{k!} i_{\la
!}\left(\sum\limits_{\square\in \lambda} c_{\square}^k t^k \right
)
=\frac{1}{k!}\sum\limits_{\square\in \lambda} c_{\square}^k
t^k[\xi_{\lambda}]. \label{eq:chern}
\end{eqnarray}

Let $k$ be a nonnegative integer. Denote
$$ \widetilde{\ch}_k^{[n]} =t^{n-k} \ch_{k,T}^{[n]}  \in \Hn.$$

We define an operator $\mathfrak G $, resp. $\mathfrak G_k$, in
$\End (\Hx)$ by sending $a \in \Hn$ to $a \star \sum_{k \ge 0}
\widetilde{\ch}_k^{[n]}$, resp. to $a \star
\widetilde{\ch}_k^{[n]}$, in $\Hn$ for each $n$. Similarly, we
define an operator $\mathfrak G_z$ by sending $a \in \Hn$ to $a
\star \sum_{k \ge 0} z^k \widetilde{\ch}_k^{[n]}$ for each $n$,
where $z$ is a variable. By definition, we have

$$\mathfrak G_z =\sum_{k \ge 0} \mathfrak G_k z^k.$$

Formula (\ref{eq:chern}) is equivalent to

\begin{eqnarray}  \label{eq:basic}
\mathfrak G_z ([\la]) =\sum_{\square\in\lambda} e^{z c_{\square}}
\cdot [\la].
\end{eqnarray}
We introduce

$$ \varsigma (z) =e^{z/2} -e^{-z/2}. $$

\begin{lemma}   \label{lem:basic}
Given a partition $\la =(\la_1, \la_2, \ldots)$, we have
$$\mathfrak G_z ([\la]) = \frac1{\varsigma (z)} \left(\sum_{i
=1}^\infty e^{z (\la_i -i+1/2)} -\frac1{\varsigma (z)}\right)
\cdot [\la].$$
\end{lemma}

\begin{proof}
The contents $c_\square$'s of $\la$ are:

$$ -i+1, -i+2, \ldots, -i+\la_i, \quad i =1,\ldots,\ell(\la).  $$
Noting that$$e^{z(-i+1)} +\ldots + e^{z(-i+\la_i)} =
\frac{e^{z(\la_i-i+1)}-e^{z(-i+1)}}{e^z-1}.$$
the formula (\ref{eq:basic}) can be written as

\begin{eqnarray*}
\mathfrak G_z ([\la]) = \frac1{\varsigma (z)} \sum_{i =1}^\infty
(e^{z (\la_i -i+1/2)} - e^{z( -i+1/2)}) \cdot [\la]
\end{eqnarray*}
where we have used $\la_i=0$ for $i>\ell(\la)$. Now the lemma
follows from the identity $\sum_{i =1}^\infty e^{z( -i+1/2)} =
1/\varsigma (z)$.
\end{proof}

\begin{remark}  \label{rem:dict}
Comparing the above with the computation in Lascoux-Thibon
\cite{LT}, Lemma~3.1, we observe that the Jucys-Murphy (JM)
elements of the symmetric group $S_n$ corresponds exactly to the
equivariant Chern roots of the rank $n$ bundle $\On$. A similar
observation in the non-equivariant setup has played an important
role in the study of the Chen-Ruan orbifold cohomology ring of the
symmetric product \cite{QW}. In particular, the first equivariant
Chern character of $\On$ corresponds to the conjugacy classes of
transpositions (compare \cite{FW}).
\end{remark}

It was showed in \cite{Vas} that there is a ring isomorphism
between $\Hn$ and the class algebra $R(S_n)$ of the symmetric
group $S_n$, which sends $[\la]$ to the Schur function $s_\la$
etc. Our results above and observation in Remark~\ref{rem:dict}
make the dictionary between the two rings more explicit. For the
convenience of the reader, we compile the following dictionary
table.

\vskip 2pc

\vbox{
 \tabskip=0pt \offinterlineskip
\def\tablerule{\noalign{\hrule}}
\halign to360pt{\strut#& \vrule# \tabskip=1.2em plus1.2em& \hfil#&
\vrule#&\hfil#\hfil& \vrule#& #\hfil& \vrule#& #\hfil&
\vrule#\tabskip=0pt \cr\tablerule
 && \omit\hidewidth
Hilbert Scheme Setup \hidewidth&& \omit\hidewidth  Symmetric Group
Setup \hidewidth&\cr \tablerule\tablerule\tablerule
&& \omit\hidewidth the ring $\Hn$ \hidewidth&& \omit\hidewidth the
class algebra $R(S_n)$ \hidewidth&\cr \tablerule
 &&fixed-point class $[\la]$ &&Schur function $s_\la$ &\cr\tablerule
 &&Heisenberg monomial $\mathfrak p_{-\la}$ && conjugacy class (c.c.) $K_\la$ &\cr\tablerule
 && Chern character $\widetilde{\ch}_1^{[n]}$ &&c.c. $K_{21^{n-2}}$ of transpositions  &\cr\tablerule
 && Chern character $\widetilde{\ch}_k^{[n]}$ &&$k$-th power-sum of JM elements  &\cr\tablerule
}} \vskip 2pc

\begin{remark}  \label{rem:gener}
By the dictionary, Theorem~5.10 in \cite{Wa} on the ring
generators of $R(S_n)$ implies that $\widetilde{\ch}_{k}^{[n]}$,
$0 \leq k <n$, form a set of ring generators of $\Hn$.
\end{remark}
\subsection{The moduli space $\mathcal M(m,n)$}

Denote by $\mathcal O_m$ the $T$-equivariant line bundle over
$X=\C^2$ associated to the $T$-character $t^m$, where $m \in \Z$.
Let $\mathcal M(m,n)$ be the moduli space which parameterizes all
rank-1 subsheaves of $\mathcal O_m$ such that the quotients are
supported at finitely many points of $X$ and have length $n$.
Given $\mathcal I \in \Xn$, then $\mathcal O_m \otimes \mathcal I$
is an element in $\mathcal M(m,n).$

As before, the study of the equivariant cohomology ring
$H^*_T(\mathcal M(m, n))$ leads to a ring $\Hn^{(m)}
\stackrel{\text{def}}{=} H^{2n}_T(\mathcal M(m, n))$ whose product
is denoted by $\star$. The natural identification $\mathcal M(m,
n) \cong \Xn$ leads to the natural identification of the rings
$\Hn^{(m)} \cong \Hn,$
which induces a bilinear form $\langle -, -\rangle^{(m)}_n$ on
$\Hn^{(m)}$ from $\Hn$. We introduce

\begin{eqnarray*}
\FFock^{(m)} &=&  \bigoplus_{n=0}^{\infty} \Hn^{(m)}, \quad m \in\Z \\
 \FFock &=& \bigoplus_{m\in\Z} \FFock^{(m)}.
\end{eqnarray*}
In particular we identify $\mathcal M(0,n) =\Xn$ and $\Hn^{(0)}
=\Hn$.

We denote by $S$ the isomorphism $S:\Hn^{(m)} \rightarrow
\Hn^{(m+1)}.$ This induces isomorphisms (which will be denoted by
$S$ again) $S: \FFock \rightarrow \FFock$ and $S: \FFock^{(m)}
\rightarrow \FFock^{(m+1)}$ for all $m\in \Z$. The bilinear form
on $\FFock$ induced from $\langle -, -\rangle^{(m)}_n$ on
$\Hn^{(m)}$ will be again denoted by $\langle-,-\rangle$.

\subsection{Operators on the fermionic Fock space}

By the standard boson-fermion correspondence \cite{MJD}, $\FFock$
can be identified with the fermionic Fock space (or equivalently,
the infinite wedge space). The operator $S$ is exactly the {\em
shift} operator on the fermionic Fock space. Given an operator
$\mathfrak f \in \End(\FFock)$, the number $\langle \mathfrak f
\rangle :=\langle |0\rangle, \mathfrak f |0\rangle \rangle$ is
called the {\em vacuum expectation} of $\mathfrak f$.

It has been well known (cf. e.g. \cite{MJD, Wa}) that the
completed infinite-rank general linear Lie algebra
$\widehat{gl}_\infty$ (whose standard basis is denoted by
$E_{i,j}$, $i,j\in\Z +1/2$) acts on the fermionic Fock space. As
will become clear below (cf. Lemma~\ref{lem:corres} and
Lemma~\ref{lem:opnew}), the study of the equivariant intersection
theory on Hilbert schemes naturally leads to the following
operator in $\End (\FFock)$:

\begin{eqnarray*}
\mathfrak H (z)  =\frac1{\varsigma (z)} \sum_{k \in \Z +\frac12}
e^{kz} E_{k,k}.
\end{eqnarray*}
The operator is further expanded as

\begin{eqnarray*}
\mathfrak H (z)  =\sum_{k=0}^\infty {\mathfrak H_k} z^k
\label{eq:taylor}.
\end{eqnarray*}

\begin{remark}
The operator $\mathfrak H (z)$ in a somewhat different form has
appeared in the study of the class algebras of the symmetric
groups and wreath products \cite{LT, Wa}, and it affords a compact
expression in terms of vertex operators.
\end{remark}

On the other hand, the following operators

\begin{eqnarray} \label{eq:epsilon}
\varepsilon_r(z) =\sum_{k \in \Z +\frac12} e^{z(k-r/2)} E_{k-r,k}
+ \frac{\delta_{r,0}}{\varsigma (z)}, \quad r\in \Z
\end{eqnarray}
have been introduced by Okounkov and Pandharipande, and they have
played a fundamental role in the study of the stationary
Gromov-Witten invariants \cite{OP}. The operators
$\varepsilon_r(z)$, $r \in\Z$, satisfy the following identities:

\begin{eqnarray}
 {[} \mathfrak p_k, \varepsilon_r(z)]
 &=&\varsigma (kz) \varepsilon_{k+r}(z)  \label{eq:commu} \\
 S^{-1} \varepsilon_0(z) S
  &=& e^z \varepsilon_0(z) \label{eq:conj} \\
 \varepsilon_0( z) ([\la])
  &=& \sum_{i=1}^\infty e^{z(\la_i -i+1/2)} \cdot [\la]
  \label{eq:eigenvalue}
\end{eqnarray}
for all partitions $\la =(\la_1, \la_2, \ldots)$ (here we recall
$[\la] \in \Hx =\FFock^{(0)}$).

\begin{lemma}  \label{lem:ident}
 We have the following identities for operators on $\FFock$:

\begin{eqnarray}
\mathfrak H (z) &=& \frac1{\varsigma (z)} \left(\varepsilon_0(z)
-\frac1{\varsigma (z)} \text{\rm I} \right) \label{eq:identify}  \\
S^{-m} \mathfrak H (z) S^m &=& e^{mz}\mathfrak H (z) +
\frac{e^{mz}-1}{\varsigma(z)^2} \text{\rm I}, \quad m\in\Z
\label{eq:op}
\end{eqnarray}
\end{lemma}

\begin{proof}
The first identity follows from the definitions. The second one
follows from the first one and (\ref{eq:conj}).
\end{proof}

\begin{lemma}   \label{lem:corres}
 As an operator on $\FFock^{(0)} =\Hx$, we have the identification
$$\mathfrak G_z  = \mathfrak H(z).$$
\end{lemma}

\begin{proof}
Follows from Lemma~\ref{lem:basic}, formula (\ref{eq:eigenvalue})
and Lemma~\ref{lem:ident}.
\end{proof}

In light of the interpretation of $\mathfrak H(z)$ in
Lemma~\ref{lem:corres} and the role of $\varepsilon_0(z)$ in
\cite{OP}, the identity (\ref{eq:identify}) defines the
Gromov-Witten/Hilbert correspondence.

%
%
%

\section{The $n$-point functions of equivariant intersection numbers}
\label{sec:Npoint}

\subsection{Equivariant intersection numbers on Hilbert schemes}

Keeping in mind Remarks~\ref{rem:same} and \ref{rem:gener}, one
way to understand the equivariant intersection theory on $\Xn$ is
to study the products

$$\widetilde{\ch}_{k_1}^{[n]} \star \widetilde{\ch}_{k_2}^{[n]} \star \ldots  \star
\widetilde{\ch}_{k_s}^{[n]}$$
for arbitrary nonnegative integers $k_1, \ldots, k_s$, $s \ge 1$.

For partitions $\la$ and $\mu$ of $n$, we may consider the
intersection numbers

\begin{eqnarray}   \label{eq:twopoint}
\langle \la,\widetilde{\ch}_{k_1}^{[n]}  \cdots
\widetilde{\ch}_{k_s}^{[n]},\mu \rangle_n \stackrel{\text{def}}{=}
\left\langle \mathfrak p_{-\la},\widetilde{\ch}_{k_1}^{[n]} \star
\widetilde{\ch}_{k_2}^{[n]} \star \ldots  \star
\widetilde{\ch}_{k_s}^{[n]}\star \mathfrak p_{-\mu}
\right\rangle_n
\end{eqnarray}

The main idea below and in later sections is to form some suitable
generating functions of these intersection numbers and then
reformulate them in terms of the operator formalism.

\subsection{The $n$-point function}

Given partitions $\la$ and $\mu$ of $n$, we organize the
intersection numbers $\langle \la,\widetilde{\ch}_{k_1} \ldots
\widetilde{\ch}_{k_N},\mu \rangle_n$ defined in
(\ref{eq:twopoint}) into a generating function by introducing the
{\em $N$-point function}

$$G_{\la, \mu} (z_1,\ldots, z_N)
= \sum_{k_1,\ldots, k_N} z_1^{k_1}\ldots z_N^{k_N}\langle
\la,\widetilde{\ch}_{k_1}^{[n]} \cdots
\widetilde{\ch}_{k_N}^{[n]},\mu \rangle_n.$$
This can be reformulated in an operator form:

\begin{eqnarray*}
G_{\la, \mu} (z_1,\ldots, z_N)
 &=& \langle \mathfrak p_{-\la},\mathfrak
G_{z_1} \dots \mathfrak G_{z_N} \mathfrak p_{-\mu} \rangle_n  \\
&=& \langle \mathfrak p_{-\la},\mathfrak H(z_1) \dots \mathfrak H
(z_N) \mathfrak p_{-\mu} \rangle_n.
\end{eqnarray*}

We also define similarly

\begin{eqnarray*}
F^\bullet_{\la, \mu} (z_1,\ldots, z_N) = \langle \mathfrak
p_{-\la},\varepsilon_0(z_1) \dots \varepsilon_0(z_N) \mathfrak
p_{-\mu} \rangle_n.
\end{eqnarray*}
According to \cite{OP}, Section~3, $F^\bullet_{\la, \mu}
(z_1,\ldots, z_N)$ has the interpretation as the $N$-point {\em
disconnected} series of stationary Gromov-Witten invariants of
$\mathbb P^1$ relative to $0,\infty \in \mathbb P^1$.

Note that, in particular for $N=0$, we have defined $G_{\la, \mu}
( )$ and $F^\bullet_{\la, \mu} ( )$. It is easy to see that

$$G_{\la, \mu} ( ) =F^\bullet_{\la, \mu} ( )
=\frac{\delta_{\la,\mu}}{{\mathfrak z}_\la}.$$
It follows from (\ref{eq:identify}) and Lemma~\ref{lem:corres}
that the 1-point function is given by

\begin{eqnarray}  \label{eq:relation}
 G_{\la, \mu} (z)
 = \frac1{\varsigma(z)}
 \left( F^\bullet_{\la, \mu} (z) -\frac1{\varsigma(z)} F^\bullet_{\la, \mu} ( )
 \right).
\end{eqnarray}
Therefore, it remains to compute $F^\bullet_{\la, \mu} (z).$ Note
that the somewhat simpler connected series $F^\circ_{\la, \mu}$
rather than the disconnected series has been computed in
\cite{OP}, Section~3.2.

We introduce some notations on partitions. Given two partitions
$\la$ and $\mu$, we denote by $\la +\mu$ the partition obtained
from combining the parts of $\la$ and $\mu$ and rearranging them
in a descending order. For two partitions $\la
=(1^{m_1}2^{m_2}\ldots)$ and $\mu=(1^{n_1}2^{n_2}\ldots)$, we say
$\la\subset \mu$ if $m_i \leq n_i$ for all $i$, and denote $\mu
-\la$ the partition $(1^{n_1-m_1}2^{n_2-m_2}\ldots)$.

Given a partition $\la =(\la_1, \ldots, \la_r)$, where $r
=\ell(\la)$, we denote

\begin{eqnarray*}
 \varsigma(\la,z) =\varsigma(\la_1 z) \ldots \varsigma (\la_r z).
\end{eqnarray*}
Given a subset $U \subset \underline{r} = \{1,\ldots,r\}$, we
denote by ${\la_U}$ the subpartition of $\la$ which consists of
the parts $\la_i$, $i\in U$.

\begin{proposition} \label{prop:onepoint}

Let $\la =(\la_1, \ldots, \la_r)$ and $\mu =(\mu_1, \ldots,
\mu_s)$ be partitions of $n$, where $r =\ell(\la)$ and
$s=\ell(\mu)$. We have

 \begin{eqnarray*}
 F^\bullet_{\la, \mu} (z)
 =\sum_{U}
 \frac{\varsigma (\la_U, z) \, \varsigma (\la_U +\mu-\la,z)}{\mathfrak z_\la
 \mathfrak z_{\la_U  +\mu-\la} \varsigma (z)}
 \end{eqnarray*}
summed over subsets $U \subset \underline{r}$ such that $\la
\subset \la_U+\mu$.
\end{proposition}

\begin{proof}
We denote by $[\mathfrak p_{\la,U}, \varepsilon_0(z)]$ the
multi-commutator
$[ \cdots [\mathfrak p_{\la_b}, [\mathfrak
p_{\la_a},\varepsilon_0(z)]]\cdots],$
if we write $U=\{a,b,\ldots\}$. Note that the multi-commutator is
independent of the ordering of elements in $U$ since the
$\mathfrak p_k$'s ($k>0$) commute with each other.

By moving $\varepsilon_0 (z)$ to the left whenever possible, we
have

\begin{eqnarray}
\mathfrak p_{\la_r} \ldots  \mathfrak p_{\la_1} \varepsilon_0(z)
  &=& \sum_{U \subset \underline{r}} [\mathfrak p_{\la,U}, \varepsilon_0(z)]
 \prod_{i\in\underline{r}\backslash U} \mathfrak p_{\la_i} \nonumber \\
 &=& \sum_{U \subset \underline{r}}  \varsigma (\la_U, z)  \varepsilon_{|\la_U|}(z)
 \prod_{i\in\underline{r}\backslash U} \mathfrak p_{\la_i}  \label{eq:commute}
 \end{eqnarray}
where we have repeatedly used (\ref{eq:commu}). It follows that
\begin{eqnarray}
 F^\bullet_{\la, \mu} (z)
  &=& \frac1{\mathfrak z_\la \mathfrak z_\mu}
  \langle \mathfrak p_{\la_r} \ldots  \mathfrak p_{\la_1} \varepsilon_0(z)
  \mathfrak p_{-\mu_1} \ldots  \mathfrak p_{-\mu_s} \rangle \nonumber \\
 &=& \frac1{\mathfrak z_\la \mathfrak z_\mu}
 \sum_{U \subset \underline{r}}\varsigma (\la_U, z)
   \left \langle \varepsilon_{|\la_U|}(z)
 \prod_{i\in\underline{r}\backslash U} \mathfrak p_{\la_i} \cdot
  \mathfrak p_{-\mu_1} \ldots  \mathfrak p_{-\mu_s} \right
  \rangle  \label{eq:step}
\end{eqnarray}

Denote by $\la^-_U$ to be the partition which consists of the
parts $\la_i$, $i \in \underline{r}\backslash U$. Apparently, the
vacuum expectation

$$\left \langle \varepsilon_{|\la_U|}(z)
 \prod_{i\in\underline{r}\backslash U} \mathfrak p_{\la_i} \cdot
  \mathfrak p_{-\mu_1} \ldots  \mathfrak p_{-\mu_s} \right \rangle
  =0 \quad\text{ unless }\la_U^- \subset \mu.$$
If $\la_U^- \subset \mu$, or equivalently if $\la \subset \la_U
+\mu$, then $\mu - \la_U^- = (\la_U  +\mu) -\la$, and we can show
by induction that

\begin{eqnarray}  \label{eq:induction}
\prod_{i\in\underline{r}\backslash U} \mathfrak p_{\la_i} \cdot
\mathfrak p_{-\mu_1} \ldots  \mathfrak p_{-\mu_s} \vac
=\frac{\mathfrak z_\mu}{\mathfrak z_{\la_U +\mu-\la}} \mathfrak
p_{-\mu_{a_1}} \ldots \mathfrak p_{-\mu_{a_t}} \vac
\end{eqnarray}
where we have denoted $\mu - \la_U^- =(\mu_{a_1},\ldots,
\mu_{a_t}).$

By (\ref{eq:step}) and (\ref{eq:induction}), we have

\begin{eqnarray}  \label{eq:step2}
 F^\bullet_{\la, \mu} (z)
  = \frac1{\mathfrak z_\la \mathfrak z_\mu}\sum_{U \subset \underline{r}} \frac{\mathfrak
  z_\mu}{\mathfrak z_{\la_U +\mu-\la}} \varsigma (\la_U, z)
   \left \langle \varepsilon_{|\la_U|}(z)\mathfrak p_{-\mu_{a_1}} \ldots
  \mathfrak p_{-\mu_{a_t}} \right \rangle
\end{eqnarray}

Similar to (\ref{eq:commute}), we now move
$\varepsilon_{|\la_U|}(z)$ to the right in  $\left \langle
\varepsilon_{|\la_U|}(z)\mathfrak p_{-\mu_{a_1}} \ldots
  \mathfrak p_{-\mu_{a_t}} \right \rangle$ whenever possible. Note that if
$\varepsilon_K(z)$ ($K>0)$ results from such a move, then the
corresponding vacuum expectation is zero. Therefore

\begin{eqnarray}  \label{eq:comm2}
\left \langle \varepsilon_{|\la_U|}(z)\mathfrak p_{-\mu_{a_1}}
\ldots \mathfrak p_{-\mu_{a_t}} \right \rangle
  &=&
 \left \langle [\cdots [\varepsilon_{|\la_U|}(z), \mathfrak
 p_{-\mu_{a_1}}], \cdots \mathfrak p_{-\mu_{a_t}} ]\right \rangle \nonumber  \\
  &=& \varsigma (\la_U  +\mu-\la,z) \langle \varepsilon_0(z) \rangle \nonumber \\
  &=& \varsigma (\la_U  +\mu-\la,z) /\varsigma (z).
\end{eqnarray}

Now the proposition follows from (\ref{eq:step2}) and
(\ref{eq:comm2}).
\end{proof}

\begin{theorem}
 The 1-point function $G_{\la, \mu} (z)$ is given by
$$G_{\la, \mu} (z)
 = \frac1{{\mathfrak z}_\la\, \varsigma(z)^2}
 \left(\sum_{U}
 \frac{\varsigma (\la_U, z) \, \varsigma (\la_U  +\mu-\la,z)}{
 \mathfrak z_{\la_U  +\mu-\la}}
  -  \delta_{\la,\mu}
 \right)
$$
where $U$ runs over the subsets of $\underline{\ell(\lambda)}$
such that $\la \subset \la_U +\mu$.
\end{theorem}

\begin{proof}
Follows from (\ref{eq:relation}) and
Proposition~\ref{prop:onepoint}.
\end{proof}

In general, by Lemma~\ref{lem:corres}, we have the following.

\begin{proposition} \label{prop:Npoint}
\begin{eqnarray*}
G_{\la, \mu} (z_1,\ldots, z_N) =
\frac1{\prod_{i=1}^N \varsigma (z_i)}  \sum_{U \subset
\underline{N} }  \frac{(-1)^{N-|U|}}{\prod_{i \in
\underline{N}\backslash U} \varsigma (z_i) } F^\bullet_{\la, \mu}
(z_U).
\end{eqnarray*}
where we have denoted $F^\bullet_{\la, \mu} (z_U) =\left \langle
\mathfrak p_{-\la},\prod_{i \in U} \varepsilon_0 (z_i) \mathfrak
p_{-\mu} \right \rangle_n.$
\end{proposition}

\begin{remark}
Note that the $N$-point connected series $F^\circ_{\la, \mu}$ of
Gromov-Witten invariants rather than the disconnected series
$F^\bullet_{\la, \mu}$ has been computed in \cite{OP},
Section~3.3. The strategy used in Proposition~\ref{prop:onepoint}
could be generalized to compute the $N$-point disconnected series
$F^\bullet_{\la, \mu} (z_1, \ldots, z_N)$, although the notations
would be a bit involved. Then the computation of $G_{\la, \mu}
(z_1,\ldots, z_N)$ follows by Proposition~\ref{prop:Npoint}.
\end{remark}

\subsection{The multi-point trace function}

The $q$-trace of an operator $\mathfrak f \in \End (\Hx)$ is
defined to be
$$\Tr_q \mathfrak f
\stackrel{\text{def}}{=} \sum_\la  {{\mathfrak z}_\la}
 \langle \mathfrak p_{-\la}, \mathfrak f (\mathfrak p_{-\la}) \rangle
q^{ | \la  |}.  $$
In particular, for the identity operator $\text{I}$ on $\Hx$, we
have $\Tr_q \text{\text{I}} = 1/(q;q)_\infty $, where we have used
the notation $(a;q)_\infty =(1-a)(1-aq)(1-aq^2)\cdots$.

Our main object here is the multi-point trace function $\Tr_q(
\mathfrak G_{z_1} \cdots \mathfrak G_{z_N} )$ of a product of the
operators $\mathfrak G_{z_j}$, which encodes information about the
intersection numbers (\ref{eq:twopoint}). Here $z_1, \ldots, z_N$
are independent variables.

We denote

\begin{eqnarray*}\Theta(z) =\Theta(z;q)
 &\stackrel{\text{def}}{=}&
\eta(q)^{-3} \sum_{m \in \Z} (-1)^m q^{\frac{(m + 1/2 )^2}2} e^{(m
+1/2)z} \\
&=& ( e^{z/2}- e^{-z/2})(qz;q)_\infty (qz^{-1};q)_\infty /
(q;q)_\infty^{2}
\end{eqnarray*}
Here $\eta(q) =q^{1/24}(q;q)_\infty$ is the Dedekind eta function,
and the last identity above uses the Jacobi triple product
identity. We further define
$$\Theta^{(k)}(z)  =\frac{d^k}{dz^k} \Theta(z), \qquad k \ge 0.
$$
We agree that $\Theta^{(k)}(z)  =0$ for $k<0$.

Given a positive integer $N$, we denote $\underline{N} =\{1,2,
\ldots, N\}$. Given a finite set $U$, we denote by $S_U$ the
symmetric group on $U$. In particular, $S_{\underline{N}} =S_N$.
Given $U=\{u_1,\ldots, u_k\} \subset \underline{N}$ with $u_1
<\ldots < u_k$ and a permutation $\sigma \in S_U$, we denote by
$M_{U,\sigma}$ the $k \times k$ matrix whose $(i,j)$-th entry is
$\Theta^{(j-i+1)} (z_{\sigma{u_1}} \cdots z_{\sigma{u_{k-j}}})/
(j-i+1)!$. We further denote by $\Theta_{U,\sigma}$ the product
$\Theta(z_{\sigma{u_1}}) \Theta(z_{\sigma{u_1}} z_{\sigma{u_2}})
\cdots \Theta(z_{\sigma{u_1}} \cdots z_{\sigma{u_k}})$.

A main result of Bloch-Okounkov (\cite{BO}, Theorem~0.5) is the
following formula for $\Tr_q ( \varepsilon_0(z_1) \cdots
\varepsilon_0(z_N) )$ in our notation:

\begin{eqnarray} \label{eq:BO}
\Tr_q ( \varepsilon_0(z_1) \cdots \varepsilon_0(z_N))
 = \frac1{(q;q)_\infty} \sum_{\sigma \in S_N}
\frac{\det M_{\underline{N},\sigma}}{
\Theta_{\underline{N},\sigma}}.
\end{eqnarray}
The formula (\ref{eq:BO}) also computes the stationary
Gromov-Witten invariants of an elliptic curve, according to
Okounkov and Pandharipande \cite{OP}.

\begin{theorem} \label{th:theta}
We have
$$\Tr_q (  \mathfrak G_{z_1} \cdots \mathfrak G_{z_N}) =
\frac1{(q;q)_\infty\prod_{i=1}^N \varsigma (z_i)}  \sum_{U \subset
\underline{N} }  \frac{(-1)^{N-|U|} \sum_{\sigma \in S_U}
\frac{\det M_{U,\sigma}}{ \Theta_{U,\sigma}} }{\prod_{i \in
\underline{N}\backslash U} \varsigma (z_i) } .$$
\end{theorem}

\begin{proof}
By Lemma~\ref{lem:ident} and Lemma~\ref{lem:corres}, we have
$$\Tr_q (\mathfrak G_{z_1} \cdots \mathfrak G_{z_N})
= \frac1{\prod_{i=1}^N  \varsigma(z_i)} \sum_{U \subset
\underline{N} } \frac{(-1)^{N-|U|}
\Tr_q ( \prod_{i \in U}\varepsilon_0(z_i))}{\prod_{i
\in\underline{N}\backslash U} \varsigma(z_i) } .
$$
Now the theorem follows by applying (\ref{eq:BO}) by replacing the
set $\underline{N}$ by $U$.
\end{proof}

\section{Equivariant intersection and $\tau$-functions of $2$-Toda hierarchies}
\label{sec:tau}
\subsection{Hilbert schemes and $\tau$-functions}

Let $t= (t_1, t_2, \ldots)$ and $s =(s_1, s_2, \ldots)$ be two
sequences of indeterminates. Define the following half vertex
operators:

$$\Gamma_\pm (t) = \exp \left( \sum_{k>0} t_k \mathfrak p_{\pm k}/k
\right).$$
Given a partition $\mu =(\mu_1, \mu_2, \ldots)$ we write $t_\mu
=t_{\mu_1} t_{\mu_2}\cdots$. Let $x=(x_1, x_2, \ldots)$ be another
sequence of indeterminates. We introduce the following generating
function for the intersection numbers $\langle
\la,\widetilde{\ch}_{k_1} \ldots \widetilde{\ch}_{k_N},\mu
\rangle_n$ defined in (\ref{eq:twopoint}):

\begin{eqnarray*}
\tau (x,t,s)
 =\sum_n \sum_{|\la| =|\mu| =n}
 t_\la s_\mu \left\langle \la, \exp \left(\sum_{k=0}^\infty x_k
 \widetilde{\ch}_{k}^{[n]} \right), \mu \right\rangle_n
\end{eqnarray*}

Note that $\Gamma_-(s) =\sum_{n \ge 0} \sum_{|\la|=n} t_\la
\mathfrak p_{-\la}$, and $\Gamma_+(t) =\Gamma_-(t)^*$. From the
definition of $\mathfrak H_{k}$ and Lemma~\ref{lem:corres}, we see
that the $\tau$-function affords an operator formulation:

\begin{eqnarray*}
\tau (x,t,s)
 =\left \langle
 \Gamma_+(t)  \exp \left(\sum_{k=0}^\infty x_k
\mathfrak H_{k} \right) \Gamma_-(s) \right\rangle .
\end{eqnarray*}

\subsection{The Chern character operators from $\mathcal M(m,n)$}

We have a universal exact sequence:

$$0 \to \mathcal J_m \to \pi_2^*{\mathcal O}_m \to \mathcal Q_m \to 0$$
where $\pi_1, \pi_2$ are the projections of $\mathcal M(m,n)
\times X$ to the two factors. Denote by $\On_m$ the
$T$-equivariant vector bundle (i.e. torsion-free sheaf) over
$\mathcal M(m,n)$ of rank $n$ given by the push-forward
$\pi_{1*}(\mathcal Q_m)$, whose fiber over a point $\xi_\la \in
\mathcal M(m,n) \cong X^{[n]}$ is given by

\begin{eqnarray} \label{eq:fiber}
\On_m |_{\xi_\la} =\mathcal O_m \bigotimes \On  |_{\xi_\la}.
\end{eqnarray}

In the same way as defining the operators $\mathfrak G_z$ and
$\mathfrak G_k (k \ge 0)$ acting on $\Hx =\FFock^{(0)}$, we can
define the operators $\mathfrak G_z^{(m)}$ and $\mathfrak
G_k^{(m)} (k \ge 0)$ acting on $\FFock^{(m)}$ using the cup
products with  $\sum_{k \ge 0} t^{n-k} \ch_k (\On_m) z^k$ and
$t^{n-k} \ch_k (\On_m)$ respectively on $\Hn^{(m)}$. For technical
reasons below, we introduce the following modification

$$\widetilde{\ch}_k(\On_m) = t^{n-k}\ch _k(\On_m) +c^{(m)}_k$$
where the constant $c^{(m)}_k$ is defined by

$$\frac{e^{mz}-1}{\varsigma(z)^2} =mz^{-1} +\sum_{k \ge 0} c^{(m)}_k
z^k/k!.$$ Equivalently, if we define

\begin{eqnarray}
\widetilde{\mathfrak G}_z^{(m)} =\mathfrak G^{(m)}_z
+\frac{e^{mz}-1}{\varsigma(z)^2} \text{I}.
\end{eqnarray}
and further write $\widetilde{\mathfrak G}_z^{(m)} =mz^{-1}
\text{I} +\sum_{k \ge 0} \widetilde{\mathfrak G}^{(m)}_k z^k$,
then

$$\widetilde{\mathfrak G}^{(m)}_k ={\mathfrak
G}^{(m)}_k+c^{(m)}_k \text{I}, \quad k \ge 0$$
and $\widetilde{\mathfrak G}^{(m)}_k$ acts on $\Hn^{(m)}$ by the
product with $\widetilde{\ch}_k(\On_m)$. When $m=0$, we have
$c^{(0)}_k =0$ and $\widetilde{\mathfrak G}^{(0)}_k = {\mathfrak
G}^{(0)}_k = {\mathfrak G}_k$ for all $k.$

Similarly we define the equivariant intersection numbers for
$\mathcal M(m,n)$, denoted by $\langle - \rangle^{(m)}_n$. When
$m=0$ it reduces to the ones defined earlier.
\subsection{The $\tau$-functions and $\mathcal M(m,n)$}

We form the following generating function of the equivariant
intersection numbers on $\mathcal M(m,n)$:

\begin{eqnarray*}
\tau (x,t,s,m)
 =\sum_n \sum_{ |\la| =|\mu| =n}
 t_\la s_\mu \left\langle \la, \exp \left(\sum_{k=0}^\infty x_k
\widetilde{\ch}_{k} (\On_m) \right), \mu \right\rangle^{(m)}_n
\end{eqnarray*}
In particular we have $\tau (x,t,s,0) =\tau (x,t,s)$.

\begin{lemma}  \label{lem:opnew}
 As an operator on $\FFock^{(m)}$, we have the identification
$$\widetilde{\mathfrak G}^{(m)}_z =\mathfrak H (z).$$
\end{lemma}

\begin{proof}
Under the identification $S^m:\Hn =\Hn^{(0)}
\stackrel{\sim}{\rightarrow} \Hn^{(m)}$, we denote by
$[\la]^{(m)}$ the image of $[\la]$. By the identification of toric
action (\ref{eq:fiber}), the same proof as in
Lemma~\ref{lem:basic} implies that

\begin{eqnarray*}  \label{eq:Gtilde}
\widetilde{\mathfrak G}^{(m)}_z ( [\la]^{(m)})
 &=&  \left(\frac{e^{mz}}{\varsigma (z)}  \left(\sum_{i
=1}^\infty e^{z (\la_i -i+1/2)} - \frac1{\varsigma (z)}  \right) +
\frac{e^{mz}-1}{\varsigma(z)^2} \right) \cdot [\la]^{(m)}.
\end{eqnarray*}

By (\ref{eq:op}) and $S^m ([\la]) =[\la]^{(m)}$, we have

\begin{eqnarray*}
\mathfrak H (z) ([\la]^{(m)})
 = \mathfrak H (z) S^m ([\la])
 = e^{mz}  S^m \mathfrak H (z) ([\la]) +\frac{e^{mz}-1}{\varsigma(z)^2}S^m
 ([\la]).
\end{eqnarray*}
This can be rewritten as

\begin{eqnarray*}\label{eq:switch}
\mathfrak H (z) ([\la]^{(m)})
 =\left(\frac{e^{mz}}{\varsigma (z)}  \left(\sum_{i
=1}^\infty e^{z (\la_i -i+1/2)} - \frac1{\varsigma (z)}  \right) +
\frac{e^{mz}-1}{\varsigma(z)^2} \right) \cdot [\la]^{(m)}
\end{eqnarray*}
by lemma~\ref{lem:basic} and Lemma~\ref{lem:corres}. This finishes
the proof.
\end{proof}

\begin{theorem}
\begin{enumerate}
\item  The function $\tau (x,t,s,m)$ can be reformulated as:

 \begin{eqnarray*}
 \tau (x,t,s,m)
  =\left \langle
  S^{-m} \Gamma_+(t)  \exp \left(\sum_{k=0}^\infty x_k
 {\mathfrak H}_{k} \right) \Gamma_-(s) S^m \right\rangle.
 \end{eqnarray*}

\item The functions $\tau (x,t,s,m)$, $m\in\Z$, satisfies the
$2$-Toda hierarchy of Ueno-Takasaki \cite{UT}. The lowest equation
among the hierarchy reads:

 \begin{eqnarray}    \label{eq:lowest}
 \frac{\partial^2}{\partial {t_1} \partial {s_1}} \ln \tau
 (t,s,x,m)
  =\frac{\tau(t,s,x,m+1)\; \tau(t,s,x,m-1)}{\tau(t,s,x,m)^2}.
 \end{eqnarray}
\end{enumerate}
\end{theorem}

\begin{proof}
Part (1) follows from the definition of $\tau (x,t,s,m)$ and
Lemma~\ref{lem:opnew}. The second part is standard since the
operator $\mathfrak H_k$ lies in $\widehat{gl}_\infty$.
\end{proof}

By setting $t_2 =t_3 =\cdots =s_2 =s_3 = \cdots =0$ and $x_2 =x_3
=\cdots =0$ in $\tau (x,t,s)$, we obtain the following generating
function of the intersection numbers on Hilbert schemes:

\begin{eqnarray*}
\sum_n (t_1s_1 e^{x_0})^n
   \left\langle \frac1{n!} \mathfrak p_{-1}^n, \exp \left( x_1
\mathfrak G_{1} \right)\cdot \frac1{n!} \mathfrak p_{-1}^n
\right\rangle_n
\end{eqnarray*}
thanks to the fact that $\mathfrak G_0 ([\la]) =|\la|\cdot [\la]$.
Setting $u=x_0 + \ln (t_1s_1)$, we denote the above generating
function by $\tau (u,x_1)$. A simple computation reduces the Toda
equation (\ref{eq:lowest}) to the following:

\begin{eqnarray*}
e^{-u} \frac{\partial^2}{\partial u^2} \ln \tau (u,x_1)
 = \frac{\tau(u+x_1, x_1)\; \tau(u-x_1, x_1)}{\tau(u, x_1)^2}.
\end{eqnarray*}
It is interesting to observe that this $\tau$ function can also be
interpreted \cite{Ok} as generating functions of certain Hurwitz
numbers.

\end{document}